\documentclass[11pt]{amsart}
\usepackage{amscd,amssymb,graphicx}
\usepackage[nohug]{diagrams}

\begin{document}

\newtheorem{thm}{Theorem}
\newtheorem{cor}[thm]{Corollary}
\newtheorem{lemma}[thm]{Lemma}

\theoremstyle{definition}
\newtheorem{quest}[thm]{Question}
\newtheorem{ex}[thm]{Example}

\newcommand{\Chat}{\widehat{\mathbb C}}
\newcommand{\Otilde}{\widetilde{\mathcal{O}}}
\newcommand{\bR}{\mathbb R}
\newcommand{\bZ}{\mathbb Z}
\newcommand{\lcm}{\mbox{lcm}}
\newcommand{\co}{\colon\thinspace}
\newcommand{\zn}{\nu}
\newcommand{\zp}{\pi}
\newcommand{\zr}{\rho}
\newcommand{\zs}{\sigma}
\newcommand{\zv}{\varphi}
\newcommand{\zJ}{\Psi}
\newcommand{\cO}{\mathcal{O}}
\newcommand{\cQ}{\mathcal{Q}}
\newcommand{\cR}{\mathcal{R}}
\newcommand{\SR}{S_{\mathcal{R}}}
\newcommand{\expm}{\varphi}
\newcommand{\subm}{\sigma_{\mathcal{R}}}

\newcommand{\nosubsections}{\renewcommand{\thethm}{\thesection.\arabic{thm}}
           \setcounter{thm}{0}}

\title{Constructing subdivision rules from rational maps}

\author{J. W. Cannon}
\address{Department of Mathematics\\ Brigham Young University\\ Provo, UT
84602\\ U.S.A.}
\email{cannon@math.byu.edu}

\author{W. J. Floyd}
\address{Department of Mathematics\\ Virginia Tech\\
Blacksburg, VA 24061\\ U.S.A.}
\email{floyd@math.vt.edu}
\urladdr{http://www.math.vt.edu/people/floyd}

\author{W. R. Parry}
\address{Department of Mathematics\\ Eastern Michigan University\\
Ypsilanti, MI 48197\\ U.S.A.}
\email{walter.parry@emich.edu}

\begin{abstract}
This paper deepens the connections between critically finite
rational maps and finite subdivision rules. The main theorem is that
if $f$ is a critically finite rational map with no periodic critical
points, then for any sufficiently large integer $n$ the iterate
$f^{\circ n}$ is the subdivision map of a finite subdivision rule.
\end{abstract}

\thanks{This work was supported in part by NSF research grants
DMS-0104030 and DMS-0203902.} 
\date\today
\maketitle

We are interested here in connections between finite subdivision
rules and rational maps. Finite subdivision rules arose out of our
attempt to resolve Cannon's Conjecture: If $G$ is a
Gromov-hyperbolic group whose space at infinity is a 2-sphere, then
$G$ has a cocompact,  properly discontinuous, isometric action on
hyperbolic 3-space. Cannon's Conjecture can be reduced (see, for
example, the Cannon-Swenson paper \cite{CS}) to a conjecture about
(combinatorial) conformality for the action of such a group $G$ on
its space at infinity, and finite subdivision rules were developed
to give models for the action of a Gromov-hyperbolic group on the
2-sphere at infinity.

There is also a connection between finite subdivision rules and
rational maps. If $\cR$ is an orientation-preserving finite
subdivision rule such that the subdivision complex $\SR$ is a
2-sphere, then the subdivision map $\subm$ is a critically finite
branched map of this 2-sphere. In joint work \cite{ratsub} with
Kenyon we consider these subdivision maps under the additional
hypotheses that $\cR$ has bounded valence (this is equivalent to its
not having periodic critical points) and mesh approaching $0$. In
\cite[Theorem 3.1]{ratsub} we show that if $\cR$ is conformal (in
the combinatorial sense) then the subdivision map $\subm$ is
equivalent to a rational map. The converse follows from
\cite[Theorem 4.7]{expi}.

In this paper we consider the problem of when a rational map $f$ can
be equivalent to the subdivision map of a finite subdivision rule.
Since a subdivision complex has only finitely many vertices, such a
rational map must be critically finite. We specialize here to the
case that $f$ has no periodic critical points. Our main theorem,
which has also been proved by Bonk and Meyer \cite{BM} when $f$ has
a hyperbolic orbifold (which includes all but some well-understood
examples), is the following:

\begin{thm}\label{thm:main}
Suppose that $f$ is a critically finite rational map with no
periodic critical points. If $n$ is a sufficiently large positive
integer, then $f^{\circ n}$ is the subdivision map of a finite
subdivision rule.
\end{thm}

In Section~\ref{sec:fsrs} we give basic definitions about finite
subdivision rules, and in Section~\ref{sec:critfin} we give basic
definitions about critically finite maps. Theorem~\ref{thm:main} is
proved in Section~\ref{sec:mainthm}. Finally, in
Section~\ref{sec:questions} we give a few questions about realizing
rational maps by subdivision maps of finite subdivision rules. Our
work, which is inspired by Sullivan's dictionary between rational
maps and Kleinian groups, is part of an effort to further tie
together Kleinian groups and rational maps. We thank Kevin Pilgrim
for suggesting this problem to us and for numerous conversations
about rational maps and finite subdivision rules.

\section{Finite subdivision rules}
\label{sec:fsrs}

In this section we give some of the basic definitions about finite
subdivision rules. For examples and more details, see \cite{fsr}.

A {\em finite subdivision rule} $\cR$ consists of i) a {\em
subdivision complex} $\SR$, ii) a subdivision $\cR(\SR)$ of $\SR$,
and iii) a {\em subdivision map} $\subm\co \cR(\SR) \to \SR$. The
subdivision complex $\SR$ is a finite CW complex which is the union
of its closed 2-cells and which has the property that for each
closed 2-cell $\tilde{t}$ in $\SR$ there is a cell structure $t$ on
a closed 2-disk such that $t$ has at least three vertices, all of
the vertices and edges of $t$ are in the boundary of the 2-disk, and
the characteristic map $\psi_t\co t \to \SR$ with image
$\widetilde{t}$ restricts to a homeomorphism on each open cell. (The
cell complex $t$ is called a {\em tile type} of $\cR$. Similarly,
one defines an {\em edge type} corresponding to each 1-cell of
$\SR$.) The subdivision map $\subm$ is a continuous cellular map
such that restriction to each open cell is a homeomorphism.

One can use finite subdivision rules to recursively subdivide
suitable 2-complexes. Suppose $\cR$ is a finite subdivision rule. A
2-complex $X$ is an $\cR$-{\em complex} if it is the union of its
closed 2-cells and there is a continuous cellular map $f\co X\to
\SR$ (called the {\em structure map}) such that the restriction of
$f$ to each open cell is a homeomorphism. If $X$ is an $\cR$-complex
with associated map $f$, then the cell structure on $\cR(\SR)$ pulls
back under $f$ to give a subdivision $\cR(X)$ of $X$. One
recursively defines subdivisions $\cR^{n}(X)$, $n > 0$, by
$\cR^{n+1}(X) = \cR(\cR^n(X))$. In particular, the subdivision
complex $\SR$ is an $\cR$-complex and the tile types are
$\cR$-complexes. The closed 2-cells of a 2-complex are called {\em
tiles}.

Let $\cR$ be a finite subdivision rule. Then $\cR$ is {\em
orientation preserving} if there is an orientation on the union of
the open 2-cells of $\SR$ such that the restriction of $\subm$ to
each open 2-cell of $\cR(\SR)$ is orientation preserving. We say
that $\cR$ has {\em bounded valence} if there is an upper bound to
the vertex valences of the complexes $\cR^n(\SR)$, $n > 0$. If $\cR$
has bounded valence, then for each $\cR$-complex $X$ with bounded
vertex valences there is an upper bound to the vertex valences of
$\cR^n(X)$, $n
> 0$. The finite subdivision rule $\cR$ has an {\em edge pairing} if
$\SR$ is a surface, and $\cR$ has {\em mesh approaching $0$} if for
each open cover $U$ of $\SR$, there is a positive integer $n$ such
that each tile of $\cR^n(\SR)$ is contained in an element of $U$.
The finite subdivision rule $\cR$ has {\em mesh approaching $0$
combinatorially} if i) every edge of $\SR$ is properly subdivided in
some subdivision $\cR^n(\SR)$ and ii) if $e_1$ and $e_2$ are
disjoint edges of a tile type $t$, then for some positive integer
$n$ no tile of $\cR^n(t)$ intersects $e_1$ and $e_2$. This condition
is easier to check than that of mesh approaching $0$, and a finite
subdivision rule with mesh approaching $0$ combinatorially is weakly
equivalent to one with mesh approaching $0$.

\section{Critically finite maps}
\label{sec:critfin}

We give here some of the basic definitions about critically finite
maps. These came up in Thurston's topological characterization of
critically finite rational maps. For further information, see the
Thurston notes \cite{Th}, the Douady-Hubbard paper \cite{DH}, the
Pilgrim paper \cite{Pi2}, or the Pilgrim monograph \cite{Pi3}.

Let $f\co S^2 \to S^2$ be an orientation-preserving branched map. A
point $x\in S^2$ is called a {\em critical point} of $f$ if $D_f(x)
> 1$, where $D_f(x)$ is the topological degree of $f$ at $x$. We
denote the set of critical points of $f$ by $\Omega_f$ and call it
the {\em critical set} of $f$. The forward orbit of the critical set
is the {\em postcritical set} $P_f = \cup_{n > 0} f^{\circ
n}(\Omega_f)$. The map $f$ is {\em critically finite} if $f$ has
degree at least two and $P_f$ is a finite set. Two critically finite
branched maps $f,g\co S^2\to S^2$ are {\em equivalent} if there is a
homeomorphism $h\co S^2 \to S^2$ such that $h(P_f) = P_g$, $\left(
h\circ f\right)\big|_{P_f} = \left( g\circ h \right)\big|_{P_f}$,
and $h\circ f$ is isotopic rel $P_f$ to $g\circ h$.

We put an orbifold structure on the domain $S^2$ of a critically
finite branched map. Given $x\in S^2$, define $\nu_f \co S^2 \to
\bZ_+\cup \{\infty\}$ by
$$\nu_f(x) = \lcm\{D_g(y)\co
g(y)=x\ \textrm{and}\ g= f^{\circ n}\ \textrm{for some}\  n\in
\bZ_+\}.$$ Let $\cO_f$ be the orbifold $(S^2,\nu_f)$. A point $x\in
\cO_f$ with $\nu_f(x)
> 1$ is called a {\em distinguished point} or {\em ramified point}.
These are the points in $P_f$. The orbifold $\cO_f$ is {\em
hyperbolic} if there is a hyperbolic structure on $\cO_f \setminus
P_f$ such that near a distinguished point $x$ the metric is that of
a hyperbolic cone (or cusp if $\nu_f(x) = \infty$) with cone angle
$2\pi/\nu_f(x)$. One can show that the orbifold $\cO_f$ is
hyperbolic if and only if the Euler characteristic $\chi(\cO_f) = 2
- \sum_{x\in P_f}\left( 1 - \frac{1}{\nu_f(x)}\right)$ is negative.

\begin{lemma}\label{lem:ramified} Let $f\co S^2\to S^2$ be an
orientation-preserving branched map.  Let $x\in S^2$, and let
$x'=f(x)$.  Then $D_f(x)\zn_f(x)$ divides $\zn_f(x')$.
\end{lemma}
  \begin{proof} Let $n\in \bZ_+$ and $y\in S^2$, and suppose that
$f^{\circ n}(y) = x$. Let $g = f^{\circ n}$ and let $g'=f^{\circ
(n+1)}$. Then $g'(y)=x'$ and
  \begin{equation*}
D_{g'}(y)=D_f(x)D_g(y).
  \end{equation*}
Since $\zn_f(x)$ is defined as the least common multiple of all
positive integers of the form $D_g(y)$ and $\zn_f(x')$ is defined as
the least common multiple of a set of positive integers which
contains all products of the form $D_f(x)D_g(y)$, it follows that
$D_f(x)\zn_f(x)$ divides $\zn_f(x')$.
\end{proof}

Now suppose that $f\co\Chat \to \Chat$ is a critically finite
rational map. If $f$ has a periodic critical point, then the Fatou
set is not empty. If $f$ doesn't have a periodic critical point,
then $J_f = \Chat$ by Sullivan's classification \cite{Su} of the
stable components of the Fatou set. If $\cO_f$ is not hyperbolic,
then it is one of the orbifolds $(\infty,\infty)$, $(2,2,\infty)$,
$(2,3,6)$, $(2,4,4)$, $(3,3,3)$, or $(2,2,2,2)$. In each case (see,
for example, \cite{DH}), the possible rational maps are well
understood.

\begin{lemma}\label{lem:mesh}
Let $f\co \Chat \to \Chat$ be a critically finite rational map with
no periodic critical points, and let $D$ be an open topological disk
in $\Chat \setminus P_f$. Then for every $\epsilon
> 0$ there exists a positive integer $N$ such that if $n$ is an
integer with $n\ge N$, then the diameter of every component of
$f^{-n}(D)$ is less than $\epsilon$ in the spherical metric.
\end{lemma}
\begin{proof}
The rational function $f$ determines an orbifold $\cO_f =
(\Chat,\nu_f)$. Let $\cO$ be the orbifold $(\Chat,\nu_f')$, where
$\nu_f'(x) = \nu_f(x)$ if $\nu_f(x)=1$ and $\nu_f'(x) = 2 \nu_f(x)$
if $\nu_f(x) > 1$. Let $\Otilde$ be the universal covering orbifold
of $\cO$ with branched cover $p\co \Otilde \to \cO$.

In this paragraph we prove that there exists a map $g\co
\widetilde{\cO}\to \widetilde{\cO}$ which lifts $f^{-1}$; in other
words, if $x\in \widetilde{\cO}$, then $f(p(g(x)))=p(x)$.  For this,
let $\cO'=\cO\setminus P_f$ and let $\cO''=\cO\setminus (P_f\cup
f^{-1}(P_f))$.  Let $\widetilde{\cO}'=p^{-1}(\cO')$, and let
$\widetilde{\cO}''=p^{-1}(\cO'')$.  Let $p'$ denote the restriction
of $p$ to $\widetilde{\cO}'$, and let $p''$ denote the restriction
of $p$ to $\widetilde{\cO}''$.  Let $f''$ denote the restriction of
$f$ to $\cO''$.  Then $p'\co \widetilde{\cO}'\to \cO '$, $p''\co
\widetilde{\cO}''\to \cO ''$, $f''\co \cO ''\to \cO '$, and
$f''\circ p''\co \widetilde{\cO}''\to \cO '$ are covering
projections.

\begin{diagram}
\widetilde{\cO}'' &  \lDasharr   &  \widetilde{\cO}' \\
\dTo_{p''} & & \dTo_{p'}\\
\cO'' & \rTo^{f''} & \cO'\\
\end{diagram}

We wish to construct a map from $\widetilde{\cO}'$ to
$\widetilde{\cO}''$ which makes this diagram commute.  After
choosing basepoints compatibly, the lifting theorem for covering
projections implies that such a map exists if and only if
$p_*'(\zp_1(\widetilde{\cO}'))\subseteq (f''\circ
p'')_*(\zp_1(\widetilde{\cO}''))$.  The fundamental groups in
question are free groups and this containment statement follows from
the following divisibility statement: if $x,x'\in \cO$ and if
$\widetilde{x},\widetilde{x}'\in \widetilde{\cO}$ such that
$f(x)=x'$, $p(\widetilde{x})=x$ and $p(\widetilde{x}')=x'$, then
$D_{f\circ p}(\widetilde{x})$ divides $D_p(\widetilde{x}')$.  But
$D_{f\circ p}(\widetilde{x})=D_f(x) D_p(\widetilde{x})$,
$$D_p(\widetilde{x})=\begin{cases}
\zn_f(x) &\text{ if }x\notin P_f\\
2\zn_f(x) &\text{ if }x\in P_f
\end{cases},$$
and $D_p(\widetilde{x}') = 2\nu_f(x')$. Combining this with
Lemma~\ref{lem:ramified} shows that there does indeed exist a map
from $\widetilde{\cO}'$ to $\widetilde{\cO}''$ which makes the
diagram commute.  This map easily extends to a map $g\co
\widetilde{\cO}\to \widetilde{\cO}$ which lifts $f^{-1}$.

Because $f$ has degree at least two and no periodic critical points,
it follows that $\Otilde$ is hyperbolic. We next prove that $g$ is
not a homeomorphism. Since $f$ has no periodic critical points,
there must be a critical point $x$ of $f$ which is not a
postcritical point. Let $\widetilde{x}$ be a lift of $x$ in
$\Otilde$, let $y = f(x)$, and let $\widetilde{y}$ be a lift of $y$
in $\Otilde$. Then $D_p(\widetilde{x}) = 1$ and $D_p(\widetilde{y})
= 2\nu_f(y)$. Since $D_{f\circ p}(\widetilde{x})= D_f(x)
D_p(\widetilde{x}) = D_f(x)$ and $D_f(x)$ divides $\nu_f(y)$ by
Lemma~\ref{lem:ramified}, we see that $D_{f\circ p}(\widetilde{x})$
strictly divides $D_{p}(\widetilde{y})= D_g(\widetilde{y})\cdot
D_{f\circ p}(\widetilde{x}) $. Hence $D_g(\widetilde{y}) > 1$ and
$g$ is not injective in a neighborhood of $\widetilde{y}$. Hence $g$
strictly reduces distances in the hyperbolic metric. Because $\cO$
is compact and $p\co\Otilde\to\cO$ is a regular branched cover, the
group of covering transformations for $p$ acts cocompactly on
$\Otilde$. This implies that $g$ reduces distances in the hyperbolic
metric $d$ uniformly by a positive constant $c < 1$; if
$x,y\in\Otilde$, then $d(g(x),g(y)) \le c d(x,y)$. So if $n$ is a
positive integer and $\widetilde{D}$ is a lift of $D$ in $\Otilde$,
then the diameter of $g^n(\widetilde{D})$ is at most $c^n$ times the
diameter of $\widetilde{D}$. The lemma follows easily.
\end{proof}

\section{The main theorem}
\label{sec:mainthm} 

In this section we prove the main theorem, Theorem~\ref{thm:main}.
We begin with a lemma that will be used in the proof.

\begin{lemma}\label{lem:triarc}
Let $X$ be a closed topological disk with a regular tiling. If
$u,v,w$ are distinct vertices of $X$, then there is an arc in the
1-skeleton of $X$ which has initial point $u$, has terminal point
$w$, and contains $v$.
\end{lemma}
\begin{proof}
We prove this by induction on the number of tiles. If $X$ has a
single tile, then it is a polygon and the result is clear. Suppose
$n$ is a positive integer and the result holds if $X$ has a regular
tiling with $n$ tiles. Let $Y$ be a closed topological disk with a
regular tiling with $n+1$ tiles, and let $u, v, w$ be distinct
vertices in $Y$. Since the tiling of $Y$ is regular, there are at
least two tiles $t_1$ and $t_2$ such that $t_i \cap \partial Y$ is
an arc. Let $t$ be one of these such that $v$ is in the closure
$Y_1$ of $Y-t$. Then $Y_1$ is a closed topological disk with a
regular tiling with $n$ tiles. By the choice of $t$, $v$ is in
$Y_1$. If $u,v,w\in Y_1$, then it follows from induction that there
is an arc in the 1-skeleton of $Y$ from $u$ to $w$ that passes
through $v$. Now suppose that $u\not\in Y_1$ but $w\in Y_1$. Let
$u_1\ne w$ be one of the two points of $t\cap Y_1 \cap \partial Y$.
By induction there is an arc in $Y_1$ from $u_1$ to $w$ which passes
through $v$, and it can easily be extended to an arc in $Y$ from $u$
to $w$ which passes through $v$. The case $u\in Y_1$ and $w\not\in
Y_1$ is symmetric, and the case that $u,w\not\in Y_1$ is similar.
This proves the lemma.
\end{proof}

\begin{proof}[Proof of Theorem~\ref{thm:main}]
Let $\delta > 0$ such that the closures of the $\delta$
neighborhoods of the postcritical points are disjoint. We begin by
choosing a simple closed curve $\alpha$ in $\Chat$ such that
$\alpha$ contains $P_f$, for each $v\in P_f$ $\alpha$ intersects
$\overline{B_{\delta}(v)}$ in an arc, and for each $v\in P_f$
$\alpha$ intersects $\overline{B_{\delta/2}(v)}$ in an arc. We next
construct a finite CW complex $S$ with underlying space $\Chat$ such
that $P_f$ is the set of vertices of $S$ and $\alpha$ is the
1-skeleton of $S$. For each $n\in \bZ_+$, let $f^{-n}(S)$ be the
finite CW complex which is the pullback of $S$ under $f^{\circ n}$.
Note that if $\alpha \subset f^{-1}(\alpha)$, then $S$ is the
subdivision complex of a finite subdivision rule with subdivision
map $f$, and hence we are done.

Now let $\epsilon > 0$ so that $\epsilon < \delta/2$, the distance
between any two distinct components of $\alpha \setminus \cup_{v\in
P_f} B_{\delta/2}(v)$ is greater than $3\epsilon$, and there is a
point in $\Chat$ whose distance from $\alpha \cup \left(\cup_{v\in
P_f}\overline{B_{\delta}(v)}\right)$ is greater than $3\epsilon$. By
Lemma~\ref{lem:mesh}, there is a positive integer $N$ such that if
$n\ge N$ then each tile of $f^{-n}(S)$ has diameter less than
$\epsilon$. Let $n\in\bZ_+$ with $n\ge N$, and let $g = f^{\circ
n}$. It is easy to see that $P_g = P_f$.  In the next paragraph we
construct a simple closed curve $\beta$ in the 1-skeleton of
$f^{-n}(S)$ such that $\beta$ is isotopic rel $P_f$ to $\alpha$.

Choose an orientation for $\alpha$. Let $v_1,\dots,v_k$ be an
enumeration of the elements of $P_g$, labeled in cyclic order as
they appear in $\alpha$, and for convenience let $v_{k+1} = v_1$.
For each integer $i$ with $1\le i \le k$, choose a subcomplex $D_i$
of $g^{-1}(S)$ such that $D_i$ is a closed topological disk,
$B_{\delta/2}(v_i) \subset D_i$, and $D_i \subset B_{\delta}(v_i)$.
(For example, one could choose $D_i$ to be the union of the closed
star of $B_{\delta/2}(v_i)$ together with those complementary
components that are contained in $B_{\delta}(v_i)$.) For each
integer $i$ with $1\le i \le k$ let $\alpha_i$ be the subarc of
$\alpha$ from $v_i$ to $v_{i+1}$ (so $\textrm{int}(\alpha_i) \cap
\textrm{int}(\alpha_j) = \emptyset$ if $i\ne j$), and let $\beta_i$
be an arc in the 1-skeleton of $g^{-1}(S)$ such that $\beta_i\cap
D_i$ is a vertex $a_i$, $\beta_i \cap D_{i+1}$ is a vertex $b_i$,
$a_i$ and $b_i$ are the endpoints of $\beta_i$, and $\beta_i \subset
B_{\epsilon}(\alpha_i)$. For convenience we define  $b_0 = b_k$. The
choice of $\epsilon$ guarantees that $b_i \ne a_{i+1}$ for $i\in
\{0,\dots,k-1\}$, and that $\beta_i \cap \beta_j = \emptyset$ if
$i\ne j$. By Lemma~\ref{lem:triarc}, for each $i$ with $1\le i \le
k$ there is an arc $\gamma_i$ in the 1-skeleton of $D_{i}$ which
joins $a_i$ and $b_{i-1}$ and passes through $v_i$. Then the union
of the arcs $\beta_i$ and $\gamma_i$, $1\le i \le k$, is the image
of a simple closed curve $\beta$ in the 1-skeleton of $g^{-1}(S)$,
and this curve is isotopic rel $P_g$ to $\alpha$.

Since $\beta$ is isotopic to $\alpha$ rel $P_g$, there is an isotopy
$H\co \Chat\times I \to \Chat$  rel $P_g$ from the identity map on
$\Chat$ to a homeomorphism $H_1\co \Chat \to \Chat$ such that
$H_1(\alpha) = \beta$. (Here we are following the convention that
$H_t(x) = H(x,t)$.) Let $\cR$ be the finite subdivision rule with
subdivision complex $\SR = \Chat$ such that the cells of $\SR$ are
the images under $H_1$ of the cells of $S$, and with subdivision map
$\subm = H_1 \circ g$. Note that if $x$ is a vertex of $\cR(\SR) =
g^{-1}(S)$, then $g(x) \in P_g$ and so $\subm(x) = H_1(g(x)) =
g(x)$.

By the choices of $n$ and $\epsilon$, every edge of $\SR$ is
properly subdivided in $\cR(\SR)$ and, if $t$ is one of the two tile
types, no tile of $\cR(t)$ intersects two disjoint edges of $t$.
Hence $\cR$ has mesh approaching $0$ combinatorially. If $\cR$ does
not have mesh approaching $0$, then by \cite[Theorem 2.3]{fsr} $\cR$
is weakly isomorphic to a finite subdivision rule $\cQ$ with mesh
approaching $0$. Hence there are cellular homeomorphisms
$\rho_1,\rho_2 \co S_{\cQ} \to \SR$ such that $\rho_2 \circ
\sigma_{\cQ} = \subm \circ \rho_1$ and $\rho_1$ and $\rho_2$ are
cellularly isotopic. Then if we replace the subdivision map $\subm$
by $\rho_1 \circ \sigma_{\cQ} \circ \rho_1^{-1} = \rho_1 \circ
\rho_2^{-1} \circ \subm$, then we get a new finite subdivision rule,
which we still call $\cR$, which has mesh approaching $0$.
Furthermore, since $\rho_1 \circ \rho_2^{-1}$ is cellularly isotopic
to the identity map, by concatenating isotopies we still have an
isotopy $H\co \Chat \times I \to \Chat$ rel $P_g$ such that $H_0$ is
the identity map and $H_1(\alpha) = \beta$. Hence we can and do
assume that $\cR$ has mesh approaching $0$. This shows that $g$ is
equivalent to the subdivision map $\subm$. To show that $g$ is a
subdivision map of a finite subdivision rule, it suffices to show
that $g$ and $\subm$ are conjugate, since then the finite
subdivision rule $\cR$ would pull back to a finite subdivision rule
with subdivision map $g$.

We use the machinery of expansion complexes, developed in
\cite{expi}, to show that $\subm$ and $g$ are conjugate. By the
choice of $n$, there is a point in $\Chat$ whose distance from
$\beta$ is greater than  $2\epsilon$. This point is within
$\epsilon$ of a vertex of $\cR(\SR)$, and the closed star in
$\cR(\SR)$ of this vertex is in the interior of one of the two tiles
of $\SR$. By \cite[Lemma 2.4]{expi}, this tile of $\SR$ is the seed
of an expansion complex $X$ for $\cR$. Let $h\co X \to \SR$ be the
structure map for $X$, and let $\expm$ be the expansion map. Then
$X$ is homeomorphic to $\bR^2$ and $\subm \circ h = h \circ \expm$.

Since the restriction of $h$ to $h^{-1}(\SR \setminus P_g)$ is a
covering map, the maps $\{H_t \circ g\co t\in I\}$ lift to a
1-parameter family $\{\Psi_t\co t\in I\}$ of maps from $X$ to $X$
with $\Psi_1 = \expm$. Let $\Psi = \Psi_0$. Then $\Psi$ is a
homeomorphism and $g \circ h = h \circ \Psi$. The standard conformal
structure on $\Chat$ gives a nonsingular partial conformal structure
on $\SR$, and so this lifts to a conformal structure on $X$.
Furthermore, $\Psi$ is conformal with respect to this lifted
conformal structure. It now follows from \cite[Theorems 6.6 and
6.7]{expi} that there is an $\cR$-complex $X'$ with structure map
$h'\co X' \to \SR$ such that $X$ and $X'$ have the same underlying
space, $h' \circ \rho = h$ where $\rho(x) = \lim_{n\to\infty}
\Psi^{-n}(\expm^n(x))$, and $\Psi$ is the expansion map for $X'$.

We prepare to show that if $x,y\in X$, then $h(x)=h(y)$ if and only
if $h(\zr(x))=h(\zr(y))$.  For each positive integer $n$ set
$\zr_n=\zJ^{-n}\circ \zv^n$.  By continuity it suffices to prove for
every positive integer $n$ that if $x,y\in X$, then $h(x)=h(y)$ if
and only if $h(\zr_n(x))=h(\zr_n(y))$.  To begin the proof of this,
let $k$ be a nonnegative integer, and let
  \begin{equation*}
S_\cR^{(k)}=\zs_\cR^{-k}(S_\cR\setminus P_g)\quad\text{ and }\quad
X^{(k)}=h^{-1}(S_\cR^{(k)}).
  \end{equation*}
The restriction of $h$ to $X^{(0)}$ is a covering projection onto
$S_\cR^{(0)}$, and the restriction of $H$ to $S_\cR^{(0)}\times I$
lifts to a homotopy $X^{(0)}\times I\to X^{(0)}$ such that this
homotopy extends to a homotopy $\widetilde{H}\co X\times I\to X$
with $\widetilde{H}_0=\text{id}$.  We now have the following
commutative diagram.
  \begin{equation*}
\begin{CD}
X^{(k)}\times I @>\zv^k\times \text{id}>> X^{(0)}\times I
@>\widetilde{H}>>
X^{(0)}@<\zv^k<< X^{(k)} \\
@VVh\times \text{id}V  @VVh\times \text{id}V  @VVhV @VVhV \\
S_\cR^{(k)}\times I  @>\zs_\cR^{k}\times \text{id}>>
S_\cR^{(0)}\times I @>H>>
S_\cR^{(0)} @<\zs_\cR^k<< S_\cR^{(k)}\\
\end{CD}
  \end{equation*}

In this paragraph we show that $\zv^{-k}\circ \widetilde{H}_1\circ
\zv^k$ is a lift of a map from $S_\cR$ to $S_\cR$.  The restriction
of $\zs_\cR^k$ to $S_\cR^{(k)}$ is a covering projection onto
$S_\cR^{(0)}$ and the map $H\circ (\zs_\cR^k\times \text{id})$ lifts
to a homotopy $S_\cR^{(k)}\times I\to S_\cR^{(k)}$ which is the
restriction of a homotopy $L\co S_\cR\times I\to S_\cR$ with
$L_0=\text{id}$.

\begin{diagram}
     &    &     &     & S_\cR \\
     &       &     &\ruTo(4,2)^L    & \dTo_{\zs_\cR^k}\\
S_\cR\times I &\rTo^{\,\ \zs_\cR^k\times \text{id}} &S_\cR\times I &\rTo^H &S_\cR\\
\end{diagram}

So both $L\circ (h\times \text{id})$ and $h\circ \zv^{-k}\circ
\widetilde{H}\circ (\zv^k\times \text{id})$ define continuous maps
from $X^{(k)}\times I$ to $S_\cR^{(k)}$ whose restrictions to
$X^{(k)}\times \{0\}$ equal $h$ and
  \begin{equation*}
\begin{aligned}
\zs_\cR^k\circ L\circ (h\times \text{id}) & =H\circ (\zs_\cR^k\times
\text{id})\circ (h\times \text{id})\\
& =h\circ \widetilde{H}\circ
(\zv^k\times \text{id})\\
 & =h\circ \zv^k\circ \zv^{-k}\circ \widetilde{H}\circ (\zv^k\times
\text{id})\\
 & =\zs_\cR^k\circ h\circ \zv^{-k}\circ
\widetilde{H}\circ (\zv^k\times \text{id}).
\end{aligned}
  \end{equation*}
Hence $L\circ (h\times \text{id})=h\circ \zv^{-k}\circ
\widetilde{H}\circ (\zv^k\times \text{id})$, and so $\zv^{-k}\circ
\widetilde{H}_1\circ \zv^k$ is a lift of $L_1$.  Thus $\zv^{-k}\circ
\widetilde{H}_1\circ \zv^k$ is a lift of a map from $S_\cR$ to
$S_\cR$.

Next note that since $H_1\co S_\cR\to S_\cR$ is a homeomorphism
which is homotopic to the identity map rel $P_g$, so is its inverse.
Let $G\co S_\cR\to S_\cR$ be a homotopy rel $P_g$ from the identity
map to $H_1^{-1}$.  The argument of the previous paragraph with $G$
instead of $H$ shows that $\zv^{-k}\circ \widetilde{G}_1\circ \zv^k$
is a lift of a map from $S_\cR$ to $S_\cR$.  Furthermore
$\widetilde{H}_1=(\widetilde{G}_1)^{-1}$ and since $g=H_1^{-1}\circ
\zs_\cR=G_1\circ \zs_\cR$, the case $k=1$ shows that
$\zJ=\widetilde{G}_1\circ \zv$ lifts $g$.

So
  \begin{equation*}
\begin{aligned}
\zr_n & =\zJ^{-n}\circ \zv^n=(\widetilde{G}_1\circ \zv)^{-n}\circ
\zv^n=(\zv^{-1}\circ (\widetilde{G}_1)^{-1})^n\circ \zv^n\\
 & =(\zv^{-1}\circ \widetilde{H}_1\circ \zv)\circ
    (\zv^{-2}\circ \widetilde{H}_1\circ \zv^2)\circ
    (\zv^{-3}\circ \widetilde{H}_1\circ \zv^3)\circ\cdots \circ
    (\zv^{-n}\circ \widetilde{H}_1\circ \zv^n).
\end{aligned}
  \end{equation*}
We have seen that every parenthesized function $F$ has the property
that if $x,y\in X$ with $h(x)=h(y)$, then $h(F(x))=h(F(y))$.  Hence
$\zr_n$ also has this property.  Similarly,
  \begin{equation*}
\begin{aligned}
\zr_n^{-1} & =\zv^{-n}\circ \zJ^{n}=\zv^{-n}\circ
(\widetilde{G}_1\circ \zv)^{n}
=\\
 & =(\zv^{-n}\circ \widetilde{G}_1\circ \zv^n)\circ\cdots \circ
    (\zv^{-3}\circ \widetilde{G}_1\circ \zv^3)\circ
    (\zv^{-2}\circ \widetilde{G}_1\circ \zv^2)\circ
    (\zv^{-1}\circ \widetilde{G}_1\circ \zv).
\end{aligned}
  \end{equation*}
As for $\zr_n$, it follows that if $x,y\in X$ with $h(x)=h(y)$, then
$h(\zr_n^{-1}(x))=h(\zr_n^{-1}(y))$.  Thus if $x,y\in X$, then
$h(x)=h(y)$ if and only if $h(\zr (x))=h(\zr(y))$.

The branched covering $h\co X\to \SR$ is a quotient map, and by the
paragraph above the map $h'=  h \circ \rho^{-1}\co X\to \SR$
descends to the quotient space to a map $\mu\co \SR \to \SR$.
Similarly, the branched covering $h'\co X\to \SR$ is a quotient map,
and the map $h\co X\to \SR$ descends to the quotient space to a map
$\mu'\co \SR \to \SR$. It is easily seen that $\mu' = \mu^{-1}$.
Note that
$$g\circ \mu \circ h' = g\circ h = h \circ \Psi = \mu \circ h' \circ
\Psi = \mu \circ \subm \circ h'.$$ Let $t_1$ and $t_2$ be the two
open tiles of $\SR$, and let $\widetilde{t_1}$ and $\widetilde{t_2}$
be components of their preimages under $h'$. Since the restrictions
of $h'$ to $\widetilde{t_1}$ and to $\widetilde{t_2}$ are
homeomorphisms, $g\circ \mu = \mu \circ \subm$ on $t_1 \cup t_2$.
Since $t_1\cup t_2$ is dense in $\SR$, by continuity $g\circ \mu =
\mu \circ \subm$ and so $g$ and $\subm$ are conjugate.
\end{proof}

\section{Questions}
\label{sec:questions} \nosubsections

We close with a few questions.

\begin{quest} Is Theorem~\ref{thm:main} true for all positive
integers $n$, instead of just for $n$ sufficiently large? While the
statements here and in \cite{BM} are for $n$ sufficiently large, we
know of no obstruction for small values of $n$. It would be very
interesting if the theorem were not true for all positive integers
$n$.
\end{quest}

\begin{quest} Is Theorem~\ref{thm:main} true without the assumption
that there are no periodic critical points? Most of the machinery
for studying finite subdivision rules, including the expansion
complex machinery used here, was developed for applications where
the assumption of bounded valence was not restrictive. We do not
know if the proof given here generalizes to the case where $f$ has
periodic critical points.
\end{quest}


\begin{thebibliography}{1234}

\bibitem{BM}
M.~Bonk and D.~Meyer, \emph{Topological rational maps and
subdivisions}, in preparation.

\bibitem{fsr}
J.~W.~Cannon, W.~J.~Floyd, and W.~R.~Parry, \emph{Finite subdivision
rules}, Conform. Geom. Dyn. \textbf{5} (2001), 153--196
(electronic).

\bibitem{ratsub}
J.~W.~Cannon, W.~J.~Floyd, R.~Kenyon, and W.~R.~Parry,
\emph{Constructing rational maps from subdivision rules}, Conform.
Geom. Dyn. \textbf{7} (2003), 76--102 (electronic).

\bibitem{expi}
J.~W.~Cannon, W.~J.~Floyd, and W.~R.~Parry, \emph{Expansion
complexes for finite subdivision rules I}, Conform. Geom. Dyn.
\textbf{10} (2006), 63--99 (electronic).

\bibitem{CS}
J.~W.~Cannon and  E.~L.~Swenson, \emph{Recognizing constant
curvature groups in dimension 3}, Trans. Amer. Math. Soc.
\textbf{350} (1998), 809--849.

\bibitem{DH}
A.~Douady and J.~H.~Hubbard, \emph{A proof of Thurston's topological
characterization of rational functions}, Acta Math. \textbf{171}
(1993), 263--297.

\bibitem{Mi}
J.~Milnor,
Dynamics in One Complex Variable: Introductory Lectures,
Vieweg, Braunschweig/Wiesbaden, 1999.

\bibitem{Pi2}
K.~M.~Pilgrim,
\emph{Canonical Thurston obstructions}, Adv. Math.
\textbf{158} (2001), 154--168.

\bibitem{Pi3}
K.~M.~Pilgrim, \emph{Combinations of Complex Dynamical Systems},
Springer Lecture Notes in Mathematics 1827, Springer-Verlag, Berlin
Heidelberg New York, 2003.

\bibitem{Su}
D.~Sullivan,
\emph{Quasiconformal homeomorphisms and dynamics I.
Solution of the Fatou-Julia problem on wandering domains}, Ann.
Math. {\bf 122} (1985), 401--418.

\bibitem{Th}
W.~P.~Thurston, \emph{Lecture notes}, CBMS Conference, University of
Minnesota at Duluth, 1983.
\end{thebibliography}
\end{document}